\documentstyle[12pt]{article}
\title{Circle actions and ${\bf Z}/k$-manifolds}
\author{Weiping Zhang\thanks{Partially
supported by  the
MOEC  and the 973 project.}\\
Nankai Institute of Mathematics\\  Nankai University\\
Tianjin 300071,\  P. R. China\\
{\it E-mail: weiping@nankai.edu.cn}}
\date{  }
\begin{document}
\maketitle

\begin{abstract}
We establish an $S^1$-equivariant index theorem for Dirac
operators on ${\bf Z}/k$-manifolds.  As an application, we
generalize the Atiyah-Hirzebruch vanishing theorem for
$S^1$-actions on closed spin manifolds to the case of
 ${\bf Z}/k$-manifolds.
\end{abstract}
\underline{R\'esum\'e fran\c{c}ais} On \'etablit un th\'eor\`eme
d'indice $S^{1}$-\'equivariant pour les op\'erateurs de Dirac sur des
${\bf Z}/k$ vari\'et\'es. On donne une application de ce
r\'esultat, qui g\'en\'eralise le th\'eor\`eme d'Atiyah-Hirzebruch sur
les actions de $S^{1}$ aux ${\bf Z}/k$ vari\'et\'es. \\

$\ $

\noindent \underline{Titre fran\c{c}ais} Actions du cercle et
${\bf Z}/k$ vari\'et\'es.

$$\ $$

{\bf \S 1. $S^1$-actions and the vanishing theorem}

$\ $

Let $X$ be a closed connected smooth spin manifold admitting  a
non-trivial circle action. A classical theorem of Atiyah and
Hirzebruch [AH] states that $\widehat{A}(X)=0$, where
$\widehat{A}(X)$ is the Hirzebruch $\widehat{A}$-genus of $X$. In
this Note we  present an extension of the above result to the
case of ${\bf Z}/k$-manifolds, which were introduced by Sullivan
in his studies of geometric topology. We recall the basic
definition for completeness (cf. [F]).

$\ $

\noindent {\bf Definition 1.1} A compact connected ${\bf
Z}/k$-manifold is a compact manifold $X$ with boundary $\partial
X$, which admits a decomposition $\partial
X=\cup_{i=1}^k(\partial X)_i$ into $k$ disjoint manifolds and $k$
diffeomorphism $\pi_i:(\partial X)_i\rightarrow Y$ to a closed
manifold $Y$.

$\ $

Let $\pi:\partial X\rightarrow Y$ be the induced map.
 In what follows, we will call an
object $\alpha$ (e.g., metrics, connections, etc.) of $X$ a ${\bf
Z }/k$-object if there will be a corresponding object $\beta$ on
$Y$ such that $\alpha|_{\partial X}=\pi^*\beta$. We make the
assumption that $X$ is ${\bf Z }/k$ oriented, ${\bf Z }/k$ spin
and is of even dimension.

Let $g^{TX}$ be a ${\bf Z}/k$ Riemannian metric of $X$ which is
of product structure near $\partial X$. Let $R^{TX}$ be the
curvature of the Levi-Civita connection associated to $g^{TX}$.
Let $E$ be a ${\bf Z}/k$ complex vector bundle over $X$. Let
$g^E$ be a ${\bf Z}/k$ Hermitian metric on $E$ which is a product
metric near $\partial X$. Let $\nabla^E$ be a ${\bf Z}/k$
connection on $E$ preserving $g^E$ such that $\nabla^E$ is of
product structure near $\partial X$. Let $R^E$ be the curvature
of $\nabla^E$. Let $D^E_+:\Gamma(S_+(TX)\otimes E)\rightarrow
\Gamma(S_-(TX)\otimes E)$ be the associated  Dirac operator on $X$
and $D^E_{+,\partial X}$ (and then $D^E_Y$) be its induced Dirac
operator on $\partial X$ (and then on $Y$). Let
$\overline{\eta}(D^E_Y)$ be the reduced $\eta$-invariant of
$D^E_Y$ in the sense of  [APS]. Then
$$\widehat{A}_{(k)}(X,E) = \int_X
{\det}^{1/2}\left({{\sqrt{-1}R^{TX}/4\pi}\over
\sinh({\sqrt{-1}R^{TX}/ 4\pi})}\right){\rm
tr}\left[e^{{\sqrt{-1}\over 2\pi
}R^E}\right]-k\overline{\eta}(D^E_Y)\ \ \ {\rm mod}\ \ k{\bf
Z}\eqno (1.1)$$ does not depend on ($g^{TX}$, $g^E$, $\nabla^E$)
and determines a topological invariant in ${\bf Z}/k{\bf Z}$ (cf.
[APS] and [F]). Moreover, Freed and Melrose [FM] have proved a
mod $k$ index theorem, giving $\widehat{A}_{(k)}(X,E)\in {\bf
Z}/k{\bf Z}$ a purely topological interpretation. When $E={\bf
C}$ is the trivial vector bundle over $X$, we usually omit the
superscript $E$.

$\ $

\noindent {\bf Theorem 1.2} {\it If $X$ admits a nontrivial ${\bf
Z}/k$ circle action preserving the orientation and the Spin
structure on $TX$, then $\widehat{A}_{(k)}(X)=0$. Moreover, the
equivariant mod $k$ index in the sense of Freed and Melrose
vanishes.}

$\ $

It turns out that the original method in [AH] is difficult to
extend to the case of manifolds with boundary  to prove Theorem
1.2. Thus we will instead make use of an extension of the method
of Witten [W]. Analytic localization techniques developed by
Bismut-Lebeau [BL, Sect. 9] and their extensions to manifolds
with boundary developed in [DZ] play important roles in our proof.

$$\ $$

{\bf \S 2. A mod $k$ localization formula  for circle actions}

$\ $

We make the assumption that the ${\bf Z}/k$ circle action on $X$
lifts to a ${\bf Z}/k$ circle action  on $E$. Without loss of
generality, we may and we will assume that this ${\bf Z}/k$ circle
action  preserves $g^{TX}$, $g^E$ and $\nabla^E$. Let
$D^E_{+,APS}:\Gamma(S_+(TX)\otimes E)\rightarrow
\Gamma(S_-(TX)\otimes E)$ be the elliptic operator obtained by
imposing the standard Atiyah-Patodi-Singer boundary condition
[APS] on $D^E_+$.

Let $H$ be the Killing vector field on $X$ generated by the $S^1$
action on $X$. Then $H|_{\partial X}\subset \partial X$ induces
 a Killing vector field $H_Y$ on $Y$.
 Let ${\cal L}_H$ denote the corresponding  Lie derivative acting  on
$\Gamma(S_\pm(TX)\otimes E)$.
Then ${\cal L}_H$ commutes with
$D^E_{+,APS}$.

For any $n\in {\bf Z}$, let $F^n_\pm$ be the eigenspaces of
$\Gamma(S_\pm(TX)\otimes E)$  with respect to the eigenvalue
$2\pi n$ of ${1\over\sqrt{-1}}{\cal L}_H$. Let
$D^E_{+,APS}(n):F^n_+\rightarrow F^n_-$ be the restriction of
$D^E_{+,APS}$ on $F^n_+$.  Then $D^E_{+,APS}(n)$ is Fredholm. We
denote its  index by ${\rm ind}\,(D^E_{+,APS}(n))\in {\bf Z}$.

 Let $X_H$ (resp. $Y_H$) be the zero set of $H$ (resp.
$H_Y$) on $X$ (resp. $Y$). Then $X_H$ is a ${\bf Z}/k$-manifold
and there is a canonical map $\pi_{X_H}:\partial X_H\rightarrow
Y_H$ induced from $\pi$. We  fix a connected component
$X_{H,\alpha}$ of $X_H$, and we omit the subscript $\alpha$ if
there is no confusion.

We identify the normal bundle to $X_H$ in $X$  to the orthogonal
complement of $TX_H$ in $TX|_{X_H}$. Then $TX|_{X_H}$ admits an
$S^1$-invariant orthogonal decomposition $TX|_{X_H}=N_{m_1}\oplus
\cdots \oplus N_{m_l}\oplus TX_H$, where each $N_\gamma$,
$\gamma\in{\bf Z}$, is a complex vector bundle on which $g\in
S^1\subset{\bf C}$ acts  by multiplication by $g^\gamma$. By
using the same  notation as in [LMZ, (1.8)], we simply write that
$TX|_{X_H}=\oplus_{v\neq 0}N_v\oplus TX_H.$ Similarly, let
$E|_{X_H}$ admits the $S^1$-invariant decomposition
$E|_{X_H}=\oplus_v E_v.$

Let $S(TX_H,(\det N)^{-1})$ be the complex spinor bundle over $X_H$ associated
to the canonically induced Spin$^c$ structure on $TX_H$. It is
a ${\bf Z}/k$ Hermitian vector bundle and carries a canonically
induced ${\bf Z}/k$ Hermitian connection.


Recall that by [AH, 2.4],
 one has $\sum_v v\dim N_v \equiv 0\ {\rm mod}\ 2{\bf Z}$. Following [LMZ, (1.15)], set
$$R(q) = q^{{1 \over 2} \sum_v |v| \dim N_v} \otimes_{v>0}
 \left( {\rm Sym}_{q^v} (N_v) \otimes \det N_v\right)
 \otimes_{v<0}
  {\rm Sym}_{q^{-v}} \left(\overline{N}_v\right)
\otimes\sum_v q^v E_v = \oplus _{n} R_n q^n,$$
$$R'(q) = q^{-{1 \over 2} \sum_v |v| \dim N_v} \otimes_{v>0}
  {\rm Sym}_{q^{-v}} \left(\overline{N}_v\right)
 \otimes_{v<0}
 \left( {\rm Sym}_{q^{v}} (N_v)\otimes \det N_v\right)
\otimes\sum_v q^v E_v = \oplus _{n} R'_n q^n.$$ Then each $R_n$
(resp. $R_n'$) is a ${\bf Z}/k$ Hermitian vector bundle over $X_H$
carrying a canonically induced  ${\bf Z}/k$ Hermitian connection.
For any $n\in {\bf Z}$, let
$$D_{X_H,+}^{R_n}:\Gamma(S_+(TX_H,(\det N)^{-1})\otimes R_n)
\rightarrow \Gamma(S_-(TX_H,(\det N)^{-1})\otimes R_n)$$ be the
canonical twisted Spin$^c$ Dirac operator on $X_H$. Let
$D_{X_H,+,APS}^{R_n}$ be the corresponding elliptic operator
associated to the Atiyah-Patodi-Singer boundary condition [APS].
We will use similar notation for $R_n'$.

$\ $

\noindent {\bf Theorem 2.1} {\it For any integer $n\in {\bf Z}$, the following identities hold,}
$${\rm ind}\,D^E_{+,APS}(n) \equiv
\sum_\alpha (-1)^{\sum_{0<v}\dim N_v } {\rm
ind}\,D_{X_{H,\alpha},+,APS}^{R_n}\ \ \ {\rm mod}\ \ k{\bf
Z},\eqno(2.1)$$
$${\rm ind}\,D^E_{+,APS}(n) \equiv
\sum_\alpha (-1)^{\sum_{v<0}\dim N_v } {\rm
ind}\,D_{X_{H,\alpha}+,APS}^{R_n'}\ \ \ {\rm mod}\ \ k{\bf
Z}.\eqno(2.2)$$

{\it Proof}. For any $T\in {\bf R}$, following  Witten [W],  let
$D^E_{T,+}: \Gamma(S_+(TX)\otimes E)\rightarrow
\Gamma(S_-(TX)\otimes E)$ be the Dirac type operator defined by
$D^E_{T,+}=D^E_+ +\sqrt{-1}Tc(H)$. Let $D^E_{T,+,APS}$ be the
corresponding elliptic operator associated to the
Atiyah-Patodi-Singer boundary condition [APS]. Clearly,
$D^E_{T,+,APS}$ also commutes with the $S^1$-action. For any
integer $n$, let $D^E_{T,+,APS}(n)$ be the restriction of
$D^E_{T,+,APS}$ on $F^n_+$. Then $D^E_{T,+,APS}(n)$ is still
Fredholm. By an easy extension of [DZ, Theorem 1.2] to the current
equivariant and ${\bf Z}/k$ situation, one sees that ${\rm
ind}(D^E_{T,+,APS}(n))\ {\rm mod}\ k{\bf Z}$ does not depend on
$T\in {\bf R}$ (compare with [TZ, Theorem 4.2]).

Let $D^E_{T,+,\partial X}:\Gamma ( (S_+(TX)\otimes E)|_{\partial
X})\rightarrow \Gamma ( (S_+(TX)\otimes E)|_{\partial X})$ be the
induced Dirac type operator of $D^E_{T,+}$ on $\partial X$. For
any integer $n$, let $D^E_{T,+,\partial X}(n):F^n_+|_{\partial
X}\rightarrow F^n_+|_{\partial X}$ be the restriction of
$D^E_{T,+,\partial X}$ on $F^n_+|_{\partial X}$. Also,  the
induced Dirac operators $D^{R_n}_{+,\partial X_H}$ and
$D^{R_n}_{Y_H}$ can be defined in the same way as in Section 1.

Let $a_n>0$ be such that ${\rm
Spec}(D_{Y_H}^{R_n})\cap[-2a_n,2a_n]\subseteq\{ 0 \} .$ By
combining the techniques in [BL, Sect. 9], [BZ, Sect. 4b)] and
[LMZ, Sect. 1.2], one can prove  the following analogue of [BZ,
Theorem 3.9], stating that there exists $T_1>0$ such that for any
$T\geq T_1$,
$$\# \{ \lambda\in {\rm Spect}(D^E_{T,+,\partial
X}(n)): -a_n\leq \lambda\leq a_n \} =\dim (\ker
D^{R_n}_{+,\partial X_H}) =k\dim (\ker D^{R_n}_{Y_H}).\eqno(2.3)$$

If $\dim (\ker D^{R_n}_{Y_H})=0$, then by (2.3), one sees that
when $T\geq T_1$, $D^E_{T,+,\partial X}(n)$ is invertible. Then
${\rm ind}(D^E_{T,+,APS}(n))$ itself does not depend on $T\geq
T_1$. Moreover, by combining the techniques in [LMZ, Sect. 1.2]
and [DZ, Sect. 3], one can further prove that there exists $T_2>
0$ such that when $T\geq T_2$,
$${\rm
ind}(D^E_{T,+,APS}(n))=\sum_\alpha (-1)^{\sum_{0<v}\dim N_v } {\rm
ind}\,D_{X_{H,\alpha},+,APS}^{R_n}\eqno(2.4)$$ (compare with [DZ,
(2.13)]). From (2.4) and the mod $k$ invariance of ${\rm
ind}(D^E_{T,+,APS}(n))$ with respect to $T\in {\bf R}$, one gets
(2.1).

In general, $\dim (\ker D^{R_n}_{Y_H})$ need not be zero, and
the eigenvalues of $D^E_{T,+,\partial X}(n)$ lying in $[-a_n,a_n]$
are not easy to control. Thus the above arguments no longer apply
directly. Instead,  we observe that $\dim (\ker
(D^{R_n}_{Y_H}-a_n) )=0$, and we  use the method in [DZ] to perturb
the Dirac type operators under consideration.

To do this, let $\varepsilon>0$ be sufficiently small so that
$g^{TX}$, $g^E$ and $\nabla^E$ are of product structure on
$[0,\varepsilon]\times
\partial X\subset X$. Let $f:X\rightarrow {\bf R}$ be an
$S^1$-invariant  smooth function such that $f\equiv 1$ on
$[0,\varepsilon /3]\times \partial X$ and $f\equiv 0$ outside of
$[0,2\varepsilon /3]\times \partial X$. Let $r$ denote the
parameter in $[0,\varepsilon]$. Let $D_{X_H,-a_n,+}^{R_n}$ be the
Dirac type operator acting on $\Gamma(S_+(TX_H,(\det
N)^{-1})\otimes R_n)$ defined by
$D_{X_H,-a_n,+}^{R_n}=D_{X_H,+}^{R_n}-a_nfc(\mbox{${\partial
\over\partial r}$}).$ Let $D_{X_H,-a_n,+,APS}^{R_n}$ be the
corresponding elliptic operator associated to the
Atiyah-Patodi-Singer boundary condition [APS]. By an easy
extension of [DZ, Theorem 1.2] (compare with [TZ, Theorem 4.2]),
we see that,
$$\sum_\alpha (-1)^{\sum_{0<v}\dim N_v }{\rm ind}\,D_{X_{H,\alpha},-a_n,+,APS}^{R_n}
\equiv \sum_\alpha (-1)^{\sum_{0<v}\dim N_v }{\rm
ind}\,D_{X_{H,\alpha},+,APS}^{R_n} \ \ \ {\rm mod}\ \ k{\bf
Z}.\eqno(2.5)$$

For any $T\in {\bf R}$, let $D^E_{T,-a_n,+}: \Gamma(S_+(TX)\otimes
E)\rightarrow \Gamma(S_-(TX)\otimes E)$ be the Dirac type
operator defined by $D^E_{T,-a_n,+}=D^E_{T,+}
-a_nfc(\mbox{${\partial \over\partial r}$})$. Let
$D^E_{T,-a_n,+,APS}$ be the corresponding elliptic operator
associated to the Atiyah-Patodi-Singer boundary condition.
Let $D^E_{T,-a_n,+,APS}(n)$ be its restriction
 on $F^n_+$. Then $D^E_{T,-a_n,+,APS}(n)$ is
still Fredholm. By another extension of [DZ, Theorem 1.2], one has
$${\rm ind}\,D^E_{T,-a_n,+,APS}(n)\equiv {\rm ind}\,D^E_{T,+,APS}(n)
\ \ \ {\rm mod}\ \ k{\bf Z}.\eqno(2.6)$$

Moreover, since $D_{Y_H}^{R_n}-a_n$, which is the induced Dirac
type operator from $D_{X_H,-a_n,+}^{R_n}$ through $\pi_{X_H}$, is
invertible, by combining the arguments in [LMZ, Sect. 1.2] with
those in [DZ, Sect. 3], one deduces that there exists $T_3>0$
such that for any $T\geq T_3$, the following analogue of (2.4)
holds,
$${\rm ind}\,D^E_{T,-a_n,+,APS}(n) =
\sum_\alpha (-1)^{\sum_{0<v}\dim N_v } {\rm
ind}\,D_{X_{H,\alpha},-a_n,+,APS}^{R_n}.\eqno(2.7)$$

{}From (2.5)-(2.7) and the mod $k$ invariance of ${\rm
ind}(D^E_{T,+,APS}(n))$ with respect to $T\in {\bf R}$, one gets
(2.1).

Similarly, by taking $T\rightarrow -\infty$, one gets (2.2).
$\Box$

$$\ $$

{\bf \S 3. Proof of Theorem 1.2}

$\ $

We apply Theorem 2.1 to the case $E={\bf C}$.

First, if $X_H=\emptyset$, by Theorem 2.1,
it is obvious that for each $n\in {\bf Z}$,
$${\rm ind}\,(D_{+,APS}(n))\equiv 0\ \ \ {\rm mod}\ \ k{\bf Z}.\eqno(3.1)$$

When $X_H\neq\emptyset$, we see that $\sum_v |v| \dim N_v
>0$ (i.e., at least one of the $N_v$'s is nonzero) on each
connected component of $X_H$. Then by (2.1) and by the definition of
the $R_n$'s, we deduce that for any integer $n\leq 0$,
 (3.1) holds.
 Similarly, by (2.2) and by the definition of the $R_n'$'s, one deduces that (3.1)  holds for any integer $n\geq 0$.

In summary, for  any $n\in {\bf Z}$, (3.1) holds.

From (1.1), (3.1), by the Atiyah-Patodi-Singer index theorem
[APS], and using the obvious fact that ${\rm
ind}(D_{+,APS})=\sum_n {\rm ind}(D_{+,APS}(n) ),$ one gets
$\widehat{A}_{(k)}(X)=0$. $\Box$

$\ $

\noindent{\it Remark 3.1} By combining Theorem 2.1 with the
arguments in [LMZ, Sects. 2-4], one should be able prove an
extension of the Witten rigidity theorem, of which a
$K$-theoretic version has been worked out in [LMZ], to ${\bf
Z}/k$-manifolds. This, together with some other consequences of
Theorem 1.2, will be carried out elsewhere.

 $$\ $$

\noindent{\bf Acknowledgements} The author would like to thank
Xiaonan Ma for helpful conversations.

$$\ $$


\noindent {\bf References}

$\ $

\noindent [AH] M. F  Atiyah  and  F. Hirzebruch, Spin manifolds
and groups actions, {\it Essays on Topology and Related Topics,
M\'emoire\'e d{\'e}di{\'e} {\`a} Georges de Rham.} (Ed. A. Haefliger
and R. Narasimhan), Springer-Verlag, (1970), 18-28.

$\ $

\noindent [APS] M. F. Atiyah, V. K. Patodi and I. M. Singer,
Spectral asymmetry and Riemannian geometry I. {\it Proc.
Cambridge Philos. Soc.} 77 (1975), 43-69.





$\ $

\noindent [BL] J.-M. Bismut  and G. Lebeau, Complex immersions
and Quillen metrics, {\em Publ. Math. IHES.} 74 (1991), 1-297.

$\ $

\noindent [BZ] J.-M. Bismut and W. Zhang, Real embeddings and eta
invariant. {\it Math. Ann.} 295 (1993), 661-684.

$\ $

\noindent [DZ] X. Dai and W. Zhang,  Real embeddings and the
Atiyah-Patodi-Singer index theorem for Dirac operators. {\it Asian
J. Math.} 4 (2000), 775-794.

$\ $

\noindent [F] D. S. Freed, ${\bf Z}/k$-manifolds and families of
Dirac operators. {\it Invent. Math.} 92 (1988), 243-254.

$\ $

\noindent [FM] D. S. Freed and R. B. Melrose, A mod $k$ index
theorem. {\it Invent. Math.} 107 (1992), 283-299.

$\ $

\noindent [LMZ] K. Liu, X. Ma and W. Zhang,  Rigidity and
vanishing theorems in $K$-theory.
 {\it Commun. Anal. Geom}. 11 (2003), 121-180.

$ \ $

\noindent [TZ] Y. Tian and W. Zhang, Quantization formula for
symplectic manifolds with boundary. {\it Geom. Funct. Anal.} 9
(1999), 596-640.

$\ $

\noindent [W] E. Witten, Fermion quantum numbers in Kaluza-Klein
theory. {\it Shelter Island II: Proceedings of the 1983 Shelter
Island Conference on Quantum Field theory and the Fundamental
Problems of Physics}. (Ed. R. Jackiw et. al.), The MIT Press,
(1985), 227-277.



\end{document}